\documentclass[a4paper, 12pt]{article} 
\newdimen\paperhight
\usepackage{amssymb}
\usepackage{amsfonts}
\usepackage[usenames]{color}

\newcommand{\hf}{\frac{1}{2}}
\newcommand{\st}{\frac{1}{16}}
\newcommand{\svt}{\frac{7}{10}}
\newcommand{\pd}{\partial} 

\newcommand{\pr}{\par \vspace{3mm}\noindent [{\bf Proof}] \qquad}
\newcommand{\prend}{\hfill \qed \par \vspace{3mm}}

\newcommand{\qed}{\quad \hbox{\rule[-2pt]{3pt}{6pt}} \par \vspace{3mm}}

\newcommand{\1}{{\bf 1}} 
\newcommand{\w}{\tilde{\omega}} 
\newcommand{\C}{\mathbb C} 
\newcommand{\Z}{\mathbb Z}

\newcommand{\N}{\mathbb N}

\newcommand{\CH}{{\cal H}}

\newcommand{\al}{\alpha}

\newcommand{\ga}{\gamma}
\newcommand{\la}{\lambda}
 
\newcommand{\tr}{{\rm tr}}
\newcommand{\wt}{{\rm wt}}
\newcommand{\Hom}{{\rm Hom}}
\newcommand{\End}{{\rm End}}
\newcommand{\Res}{{\rm Res}}
\newtheorem{thm}{Theorem}[section]
\newtheorem{prn}[thm]{Proposition}
\newtheorem{dfn}[thm]{Definition}
\newtheorem{lmm}[thm]{Lemma}
\newtheorem{cry}[thm]{Corollary}

\begin{document}

\title{Intertwining operators and modular invariance}
\author{Masahiko Miyamoto
\footnote{Supported by the Grants-in-Aids for Scientific Research, 
No. 09440004 and No. 12874001, The Ministry of Education, Science and 
Culture, Japan.}} 
\date{\begin{tabular}{c}
Institute of Mathematics \cr
University of Tsukuba \cr
Tsukuba 305, Japan \cr
\end{tabular} }
\maketitle

{\large Dedicated to Professor Hiroyoshi Yamaki on his 60th birthday}\\

\begin{abstract}
We extend Zhu's theory to the case of 
intertwining operators of vertex operator algebra $V$.
Namely, we show that the space of trace functions $S^I(u, \tau)$ of 
intertwining 
operators $I$ of type ${W \choose U\quad W}$ satisfies 
modular invariance for each $U$ and $u\in U$ and we 
construct modular forms of vector type of rational weights. 
As an application, we calculate trace functions of some intertwining 
operators explicitly. 
\end{abstract}

\section{Introduction}
For a rational vertex operator algebra $V$ with central charge $c$ and 
the set of irreducible $V$-modules $ \{W^1, ..., W^m\}$, 
Zhu's theory insists that the set of trace functions 
$$ S_{W^i}(v, \tau)=z^{\wt(v)}q^{-c/24}\tr_{| W^i}Y^{W^i}(v, z)q^{L(0)}
\qquad (q = e^{2\pi i \tau})$$
for $v \in V$ satisfy some modular 
invariance ($SL(2, \Z)$-invariance) if $V$ satisfies condition 
$C_2$ (see Def. 2.7), where $Y^{W^i}(v, z)$ is the module 
vertex operator of $v$ on $ W^i $.  
Especially, if $v \in V_{[n]}$, then $SL(2, \Z)$ acts on an 
$m$-dimensional vector space 
$ \C S_{W^1}(v, \tau)+ \cdots + \C S_{W^m}(v, \tau)$ 
and $\left(S_{W^1}(v, \tau), ..., S_{W^m}(v, \tau)\right)$ 
become modular forms of 
vector type of integer weight $n$. 
This theory was extended by Dong, Li and Mason in \cite{DLiM} 
to the orbifold model, 
where $V$ has an automorphism $g$ of finite order and one consider 
the trace function $ \tr_{W^i}g q^{L(0)-c/24}$. The author has also extended 
Zhu's theory to the trace functions in many variables, 
\cite{Mi1}, \cite{Mi2}. 
As applications of these theories, 
we can construct a lot of modular forms of integer weights from 
holomorphic vertex operator algebras.  
Recently, attention has come to be paid to modular forms 
of rational weights, see \cite{BKMS} and \cite{Ib}.   
In this paper, we will show a new construction of modular 
forms of rational weights by using intertwining operators 
of vertex operator algebras. For example, 
we will construct modular forms (with a linear 
character) of weights $\hf, \frac{1}{10}, \frac{2}{5}, \frac{1}{7}$.
Actually, our proof covers 
the real weights, but we don't know such a case. 

An incentive of this research is Dedekind's $ \eta $-function
$$   \eta(\tau)=q^{1/24}\prod_{n=1}(1-q^n),  $$
which is a modular form of weight $ \hf $. 
It follows from Spinor construction \cite{FRW} that this function is 
given as a trace function 
$$   u^{-\wt(u)}q^{-\hf \frac{1}{24}}\tr_{|L(\hf, \st)} I(u, z)q^{L(0)}   $$
for some intertwining operator
$I(\ast, z)\in I\pmatrix{L(\hf, \st)\cr L(\hf, \hf)\quad L(\hf, \st)}$ 
of $L(\hf, 0)$-modules
and some element $u \in L(\hf, \hf)$ of weight $ \hf $, where 
$L(\hf, 0)$ is 2-dimensional Ising model (a rational 
Virasoro vertex operator algebra  
with central charge $ \hf $) and has three modules $L(\hf, 0)$, 
$L(\hf, \hf)$ and $L(\hf, \st)$, where the first entry denotes central charge 
and the second denotes the lowest weights.  
Since $L(\hf, \st)$ is the only one irreducible $L(\hf, 0)$-module $W$ 
satisfying $0\not=I\pmatrix{W \cr L(\hf, \hf)\quad W}$ and we have 
$ \dim I\pmatrix{L(\hf, \st)\cr L(\hf, \hf)\quad L(\hf, \st)}=1$, 
the space of trace functions 
$$ <q^{-c/24}\tr_{|W}I(u, z)q^{L(0)}: I \in I{W \choose L(\hf, \hf)\quad W}, 
W \mbox{ is an $L(\hf, 0)$-module}> $$
has dimension one for each $u \in L(\hf, \hf)$. 
This fact suggests the possibility of the extension of Zhu's theory 
to the trace functions of the intertwining operators. 
Our main purpose in this paper is to show that this is true 
and to prove that the space of trace functions $S^I(u, \tau)$ given by 
intertwining operators $I$ of type ${W \choose U\quad W}$ for each $U$ 
satisfies a modular invariance if $U$ satisfies 
a weaker condition $C_{[2,0]}$. See Theorem 4.15 and Theorem 5.1. 
Using the result (Corollary 2.13 in \cite{Li2}) given by Li, 
the proofs of these theorem are essentially the same as in \cite{Zh}. 
As we however are interested in the extension of Zhu's theory and 
the mechanics of modular invariance of trace functions of 
vertex operator algebras, 
we will pick up and repeat the necessary parts 
with the suitable modifications. \\

The author wishes to thank E. Bannai and T. Ibukiyama for their 
helpful advices. 

\section{Preliminary results}
\subsection{Vertex operator algebras}
\begin{dfn}
A vertex operator algebra is a $ \Z$-graded complex 
vector space:
$$   V= \coprod_{n \in \Z}V_n $$
satisfying $ \dim V_n<\infty $ for all $n$ and $V_n=0$ 
for $n\!<\!<\!0$.  If $ v \in V_n $ 
we write $ \wt(v)=n$ and say that $v$ is homogeneous and has (conformal) 
weight $n$. 
For each $ v \in V$ there are linear operators $ v_n \in \End(V)$, 
$ n \in \Z$ which are 
assemble into a vertex operator
$$  Y(v, z)=\sum_{n \in \Z}v(n)z^{-n-1} \in (\End V)[[z, z^{-1}]]. $$
Various axioms are imposed: \\
(1) For $ u, v \in V$, $u(n)v=0$ for $n$ sufficiently large. \\
(2) There is a distinguished {\it vacuum element} $ \1 \in V_0$ 
satisfying $Y(\1, z)=1$ 
and $Y(v, z){\1}=v+ \sum_{n \geq 2}v(-n)\1 z^{n-1}$.  \\
(3) There is a distinguished Virasoro element $ \omega \in V_2$ with generating 
function $Y(\omega, z)=\sum_{n \in \Z}L(n)z^{-n-2}$ such that the 
component operators generate a copy of the Virasoro algebra represented 
on $V$ with central charge $c$. That is 
$$  [L(m), L(n)]=(m-n)L(m + n)+ \frac{1}{12}(m^3-m)\delta_{m+n, 0}c . $$
Moreover we have  $ V_n = \{v \in V|L(0)v=nv \}$ and 
$ \frac{d}{dz}Y(v, z)=Y(L(-1)v, z)$. \\
(4)  "Commutativity" holds \\
$$  (z-w)^N[Y(v, z), Y(u, w)]=0 \mbox{ for }N\!>\!>\!0  $$
\end{dfn}

Such a vertex operator algebra may be denoted by 
the 4-tuple $(V, Y, \1, \omega)$ or more usually, by $V$. \\

It is well known that vertex operators satisfies "Associativity". 
$$(a(m)v)(r)=\sum_{i=0}^{\infty}(-1)^i{m\choose i}a(m-i)v(r+i)
-(-1)^{m+i}{m\choose i}v(m+r-i)a(i) \eqno{(2.1)}$$ 
for $a,v\in V$. 

\begin{dfn}
A module for $(V, Y, {\bf 1}, \w)$ is a 
${\bf Z}$-graded vector space $ M= \oplus_{n \geq 0}M(n)$ 
with finite dimensional homogeneous spaces $M(n)$; equipped with a formal 
power series 
$$Y^M(v, z)=\sum_{n \in {\bf Z}} v^M(n)z^{-n-1} 
\in ({\rm End}(M))[[z, z^{-1}]] $$ 
called the module vertex operator of $v$ for $v \in V$ satisfying: \\
(1)  $Y^M({\bf 1}, z)=1_M$; \\
(2)  $Y^M(\omega, z)=\sum L^M(n)z^{-n-1}$ satisfies: \\
$ \mbox{\ }$ \quad(2.a)  the Virasoro algebra relations: 
$$ [L^M(n), L^M(m)]=(n-m)L^M(n+m)+ \delta_{n+m, 0}\frac{n^3-n}{12}c, $$
$ \mbox{\ }$ \quad(2.b)  the $L(-1)$-derivative property:
$$  Y^M(L(-1)v, z)={d\over dz}Y^M(v, z), and   $$
$ \mbox{\ }$ \quad(2.c)  $L^M(0)_{M(n)}=(k_n)1_{M(n)} \mbox{ for some } 
k_n \in {\bf C}$. \\
(3)  "Commutativity" holds; 
$$  (z-w)^N[Y^M(v, z), Y^M(u, w)]=0 \mbox{ for }N\!>\!>\!0  $$
(4)  "Associativity" holds;  
$$  Y^M(u_nv, z)=Res_{w}\left((w-z)^nY^M(u, w)Y^M(v, z)-
(-z+w)^nY^M(v, z)Y^M(u, w)\right), $$
where $(-z+w)^n=\sum_{i=0}^{\infty}{n\choose i}(-z)^{n-i}w^i$ and 
${n\choose i}=\frac{n(n-1)\cdots (n-i+1)}{i!}$. 
\end{dfn}

It follows from the definitions of modules that 
if $W$ is an irreducible $V$-module then $W$ has 
a lowest weight $r$ such that $W=\oplus_{n=0}^{\infty}W_{r+n}$, where 
$L(0)$ acts on $W_{r+n}$ as a scalar $r+n$.

\begin{dfn}  A vertex operator algebra $(V, Y, \1, \omega)$ 
is called "rational" if 
it has only finitely many irreducible modules and all modules are 
completely reducible. A vertex operator algebra with a unique 
irreducible module is called "holomorphic". 
\end{dfn}

Throughout this paper, $V=\oplus_{n=0}^{\infty}V_n$ 
is a rational vertex operator algebra $(V, Y, \1, \omega)$ 
with central charge $c$ and 
$U$ is an irreducible $V$-module. We assume that $U$ 
is spanned by elements of the forms 
$L(-n_1)\cdots L(-n_t)u$ $(n_1, ..., n_t>0)$ for $u$ 
satisfying $L(n)u=0$ $(n>0)$.

Zhu has introduced the second vertex operator algebra 
$(V, Y[, ], \1, \tilde{\omega})$ associated to $V$ in 
Theorem 4.2.1 of \cite{Zh}. 

\begin{dfn}
The vertex operator $Y[a, z]$ are defined for homogeneous $a$ via the 
equality
$$  Y[a, z]=Y(a, e^{z}-1)e^{z\wt(a)} \in \End(V)[[z, z^{-1}]]  $$
and Virasoro element $ \tilde{\omega}$ is define to be $ \omega-c/24$. \\
\end{dfn}

For the proof of this fact, see \cite{Zh} or the proof of Theorem 2.1. 
We should note that Zhu has used $Y[a, z]_{2\pi i}
=Y(a, e^{2\pi iz}-1)e^{2\pi iz\wt(a)}$ 
and $ \tilde{\omega}_{2\pi i}=(2\pi i)^2(\omega-c/24)$ and Dong, Li and 
Mason have used 
$Y[a, z]_1$ in \cite{DLiM}. We adopt 
$Y[a, z]_1=Y(a, e^z)e^{\wt(a)z}$ and $ \tilde{\omega}=\omega-c/24$ because 
we can define it for a vertex operator algebra over the rational number 
field. 
It follows from the direct calculation that the differences are given by 
$$a[n]_{2\pi i}u=(2\pi i)^{-n-1}a[n]_1u.  \eqno{(2.2)} $$

Using a change of variable we calculate that 
$$\begin{array}{rl}
  v[m]&=\Res_zY[v, z]z^m \cr
      &=\Res_zY[v, \log(1+z)](\log(1+z))^m(1+z)^{-1} \cr
      &=\Res_zY(v, z)(\log(1+z))^m(1+z)^{\wt(v)-1}.  
\end{array} $$

In this paper, $ \Res_z \in \Hom(V\{z\}, V)$ is given by 
$ \Res_z(\sum_{m \in \C} a_mz^{-m-1})=a_0$. 
In particular, for $v\in V$ and $m\in \Z$, 
there are $a_i\in \C$ such that 
$$v[m]=v(m)+\sum_{i=1}^{\infty} a_iv(m+i).  \eqno{(2.3)}$$
For example we have
$$  v[0]=\sum_{i=0}^{\infty}{\wt(v)-1\choose i}v(i).  \eqno{(2.4)} $$
We also write $Y[\tilde{\omega}, z]=\sum_{n \in \Z}L[n]z^{-n-2}$ 
and set $V_{[n]}=\{v \in V| L[0]v=nv\}$. For $v \in V_{[n]}$, we denote it 
by $[\wt](v)=n$. For example, one has 
$$ L[0]=L(0)+ \sum_{n=1}^{\infty}\frac{(-1)^{n-1}}{n(n+1)}L(n) \eqno{(2.5)}$$ 
and so $ \oplus_{n\leq N}W_{r+n}=\oplus_{n\leq N}W_{[r+n]}$ for 
any irreducible $V$-module $W=\oplus_{n=0}^{\infty}W_{r+n}$.

The main operators in this paper are not vertex operators of 
elements of $V$, but 
intertwining operator of elements of $V$-module $U$. 

\begin{dfn}  Let $(W^i, Y^i)$ $(i=1, 2, 3)$ be $V$-modules. 
An intertwining operator of type $ \pmatrix{W^3\cr W^1\quad W^2}$ is a 
linear map 
$$\begin{array}{ccl}
W^1&\to &\Hom(W^2, W^3)\{z\} \cr
w_1 &\to &I(w_1, z)=\sum_{r \in \C}w_1(r)z^{-r-1} 
\end{array} $$
such that \\
(1) for $w_j \in W^j$ and $r \in \C$, $w_1(r+n)w_2=0$ for $n$ 
sufficiently large; \\
(2)  $ \frac{d}{dz}I(w_1, z)=I(L(-1)w_1, z) ;$ \\
(3)  "Commutativity" holds; 
$$  (z-w)^N\{Y^3(v, z)I(u, w)-I(u, w)Y^2(v, z)\}=0 \mbox{ for }N\!>\!>\!0 ;$$
(4)  "Associativity" holds;  
$$  I(u_1(n)v, z)=Res_{w}\left\{(w-z)^nY^3(u, w)I(v, z)-
(-z+w)^nI(v, z)Y^2(u, w)\right\}. $$
\end{dfn}

It is known (c.f. \cite{Li1}) that (3) and (4) are possible to be replaced by 
the following Jacobi identity for the operators 
$$ \begin{array}{l}
z_0^{-1}\delta(\frac{z_1-z_2}{z_0})Y^3(v, z_1)I(w_1, z_2)w_2
-z_0^{-1}\delta(\frac{z_2-z_1}{-z_0})I(w_1, z_2)Y^2(v, z_1)w_2 \cr
=z_2^{-1}\delta(\frac{z_1-z_0}{z_2})I(Y^1(v, z_0)w^1, z_2)w_2 
\end{array} \eqno{(2.6)} $$

As in the case of module actions, for $v \in V$ and $u \in W^1$ 
we have the standard consequence from (2.6): 
$$\begin{array}{l}
[Y(v, z_1), I(u, z_2)]=\Res_{z_0}z_2^{-1}\delta(\frac{z_1-z_0}{z_2})
I(Y(v, z_0)u, z_2) \cr
=I(Y(v, z_1-z_2)-Y(v, -z_2+z_1))u, z_2).  \end{array} \eqno{(2.7)}$$
Here $[Y(v, z_1), I(u, z_2)]$ denotes 
$Y^3(v, z_1)I(u, z_2)-I(u, z_2)Y^2(v, z_1)$.

By calculating the coefficients of $z_1^{-m-1}z_2^{-k-1}$ in (2.7), 
we have the following commutator formula as in the case of VOAs. 
$$ [v(m), u^I(k)]
=\sum_{j=0}^{\infty}{m\choose j}(v(j)u)^I(m+k-j), \eqno{(2.8)} $$
where $I(u, z)=\sum_{r \in \C}u^I(r)z^{-r-1}$, 
$Y^i(v, z)=\sum_{i \in \Z} v^i(m)z^{-m-1}$ 
for $i=1, 2, 3$ and \\
$[v(m), u^I(k)]$ denotes $v^3(m)u^I(k)-u^I(k)v^2(m)$.
In particular, 
$$[L(0), w^I(k)]=(\omega^1(0)u)^I(k+1)+(\omega^1(1)u)^I(k)
=(\wt(u)-k-1)w^I(k) \eqno{(2.9)}$$
and so $u^I(\wt(u)-1)$ is a grade-preserving operator 
and denoted by $o^I(u)$.\\

If $(W, Y^W)$ is a $(V, Y(, z), \1, \omega)$-module, then by 
the same arguments as in \cite{Zh}, one shows that 
$W$ become a $(V, Y[, ], \1, \tilde{\omega})$-module by 
module vertex operator $Y^W[v, z]=Y^W(v, e^z-1)e^{\wt(v)z}$. We will 
extend it to intertwining operators. 

\begin{thm}  Let $I(, z) \in I{W^3\choose W^1\quad W^2}$. 
Define the second intertwining operator $I[\ast, z]$ by 
$$  I[u, z]=I(u, e^{z}-1)e^{z\wt(u)}  \eqno{(2.10)} $$
for homogeneous $u \in W^1$. Then $I[\ast, z]$ is an 
intertwining operator of type ${W^3\choose W^1\qquad W^2}$, 
where $W^i$ are $(V, Y[, z], \1, \w)$-modules  
$(W^1, Y^1[, z]), (W^2, Y^2[, z]), (W^3, Y^3[, z])$.
\end{thm}

\pr 
As mentioned at (4.2.3) in \cite{Zh}, 
$$ a^i[m]=\Res_z(Y^i(a, z)(ln(1+z))^m(1+z)^{\wt(a)-1} \eqno{(2.11)}$$ 
for $a \in V$ and 
$i=1, 2, 3$, where $Y^i[a, z]=\sum_{m \in \Z} a^i[m]z^{-m-1}$.\\  
(1) Commutativity of $I[, z]$: 
Set $f(x, z)=(e^z-e^x)/(z-x)$. By Commutativity of 
$I(, z)$, 
$$\begin{array}{l}
0=(e^x-1-e^z+1)^N(Y^3(a, e^x-1)e^{\wt(a)x}I(u, e^z-1)e^{\wt(u)z}\cr
-I(u, e^z-1)e^{\wt(u)z}Y^2(a, e^x-1)e^{\wt(a)x}) \cr
=(x-z)^Nf(x, z)^N(Y^3[a, x]I[u, z]-I[u, z]Y^2[a, x]) 
\end{array} $$
Since $(e^z-e^x)=(z-x)(1+ \hf(x+z)+...)$, $f(x, z)$ has an inverse in 
$ \C[[x, z]]$. Hence 
$$0=(x-z)^N(Y^3[a, x]I[u, z]-I[u, z]Y^2[a, x]) \eqno{(2.12)} $$
(2) For the proofs of Associativity of $I[, z]$: 
$$(a[m]u)[r]=\sum_{i=0}^{\infty}(-1)^i{m\choose i}a[m-i]u[r+i]
-(-1)^{m+i}{m\choose i}u[m+r-i]a[i] \eqno{(2.13)}$$ 
and $L[-1]$-derivative property: 
$$I[L[-1]u, z]]=\frac{d}{dz}I[u, z], \eqno{(2.14)}$$ 
see the proof of Theorem 4.2.1 in \cite{Zh} with suitable modifications. 
\prend

We next recall the so-called condition $C_2$ introduced by Zhu \cite{Zh} and 
give a weaker condition.  

\begin{dfn}
For a $V$-module $U$, set $C_2(U)=<a(-2)u: a\in V, u\in U>$ and 
$C_{[2,0]}(U)=<a[-2]u, a[0]u: a \in V, u \in U>$.  
We call that $U$ satisfies condition $C_2$ if 
$\dim(U/C_{2}(U))<\infty$ and $U$ satisfies condition $C_{[2,0]}$ if
$\dim(U/C_{[2,0]}(U))<\infty$. 
\end{dfn}

We note $a(-n)u\in C_2(U)$ for any $n\geq 2$ since $(m-1)v(-m)
=(L(-1)v)(-m+1)$. 
From Associativity $(2.1)$, we easily have:

\begin{lmm}  $V/C_{2}(V)$ is an associative algebre 
with a product given by $v\times w=v(-1)w$. 
Also $U/C_{[2]}(U)$ is a $V/C_{[2]}(V)$-module 
whose action is given by $v(-1)u+C_2(U)$ for $v\in V, u\in U$. 
\end{lmm}

We also have the following from $(2.3)$.

\begin{lmm} If we set $C_{[2]}(U)=<a[-2]u: a\in V, u\in U>$, then 
$\dim U/C_{[2]}(U)<\infty$ if and only if $\dim U/C_2(U)$.  
\end{lmm}

\pr

\subsection{Zhu-algebra}
Let $U$ be a $V$-module. Following \cite{FZ} we define left and right 
actions of $V$ on $U$ as follows:
$$a\cdot u=\Res_x\frac{(1+x)^{\wt(a)}}{x}Y(a, x)u,  \eqno{(2.15)} $$
$$u* a=\Res_x\frac{(1+x)^{\wt(a)-1}}{x}Y(a, x)u,  \eqno{(2.16)} $$
for any homogeneous vector $a \in V$ and for any $u \in U$.  
Let $O(U)$ be the subspace of $U$ linearly spanned by all elements, 
$$ \Res_x\frac{(1+x)^{\wt(a)}}{x^2}Y(a, x)u, \eqno{(2.17)}  $$
for any homogeneous $a \in V$, $u \in U$.  Set $A(U)=U/O(U)$.  
Then it is known (Theorem 1.5.1 \cite{FZ}) that 
$A(V)$ is an associative algebra with a product $ \cdot$ and 
$A(U)$ is an $A(V)$-bimodule under the defined left and right action.
Zhu showed that $ \omega+O(V)$ is in the center of $A(V)$ and 
$A(V)$ is a finite dimensional semisimple algebra if $V$ is rational. 
As mentioned in Remark 2.9 in \cite{Li2}, 
$$ a\cdot u-u*a=\Res_x(1+x)^{\wt(a)-1}Y(a, x)u=\sum_{i=0}^{\infty}{\wt(a)-1
\choose i}a_iu=a[0]u.   \eqno{(2.18)} $$

Li recently proved the following theorem (Corollary 2.13 in \cite{Li2}), 
which was mentioned in \cite{FZ}. 

\begin{thm} If $V$ is rational, then 
there is a natural linear isomorphism 
$$\pi: I\pmatrix{W^3\cr U\qquad W^2} 
\rightarrow \Hom_{A(V)}(A(U)\otimes_{A(V)}W^2(0), W^3(0)) \eqno{(2.19)}$$
for irreducible $V$-modules $U, W^2, W^3$, where 
$W^2(0)$ and $W^3(0)$ are the top modules of $W^2$ and $W^3$, respectively. 
Here $ \pi(I)$ is given by 
$$ \pi(I)(u\otimes m)=u(\wt(u)-1+r_2-r_3)m \eqno{(2.20)} $$
for $m \in M^2(0)$, $u \in U$, $I(u, z)=\sum_{k \in \C}u(k)z^{-k-1}$ and 
$r_2, r_3$ are the lowest weights of $W^2$ and $W^3$, respectively.  
\end{thm}

\subsection{Elliptic functions}
In this section we will quote several results from \cite{Zh}. 
The Eisenstein series $G_{2k}(\tau)$ $(k=1, 2, ...)$ are series 
$$G_{2k}(\tau)
=\sum_{(m, n)\not=(0, 0)}\frac{1}{(m\tau+n)^{2k}}  \mbox{ for }k\geq 2, 
\eqno{(2.21)}  $$ 
and  
$$G_2(\tau)=\frac{\pi^2}{3}
+ \sum_{m \in \Z-\{0\}}\sum_{n \in \Z}\frac{1}{(m\tau+n)^2} 
\mbox{ for }k=1.  \eqno{(2.22)}$$
They have the $q$-expansions 
$$ G_{2k}(\tau)=2\xi(2k)+ \frac{2(2\pi\sqrt{-1})^{2k}}{(2k-1)!}
\sum_{n=1}^{\infty}\frac{n^{2k-1}q^n}{1-q^n} \eqno{(2.23)}$$ 
where $ \xi(2k)=\sum_{n=1}^{\infty}\frac{1}{n^{2k}}$ and 
$q=e^{2\pi \sqrt{-1}\tau}$.

We make use of the following normalized Eisenstein series:
$$  E_k(\tau)=\frac{1}{(2\pi\sqrt{-1})^k}G_k(\tau) \mbox{ for } k\geq 2. 
\eqno{(2.24)} $$

It is clear from (2.21) and (2.22) that 
$$ E_2({a\tau+b\over c\tau+d})=(c\tau+d)^2E_2(\tau)-\frac{c(c\tau+d)}{2\pi
\sqrt{-1}} \eqno{(2.25)}$$
for $\pmatrix{a&b\cr c&d} \in SL(2, \Z)$ and 
$E_{2k}(\tau)$ is a modular form of weight $2k$ for $k>1$. 

Set 
$$\wp_1(z, \tau)={1\over z}+
\sum_{w \in \Z \tau+\Z -\{0\}}({1\over z-w}+{1\over w}+{z\over z^w}), 
\eqno{(2.26)}$$ 
$$\wp_2(z, \tau)={1\over z^2}+ 
\sum_{w \in \Z\tau+\Z-\{0\}}({1\over (z-w)^2}-{1\over w^2}), 
\qquad \mbox{Weierstrass -$\wp$-function}  
\eqno{(2.27)}$$
and 
$$ \wp_{k+1}(z, \tau)=
-{1\over k}{d\over dz}\wp_k(z, \tau) \mbox{ for }k\geq 2.  
  \eqno{(2.28)}$$

$\frac{\wp_k(\frac{z}{2\pi\sqrt{-1}}, \tau)}{(2\pi\sqrt{-1})^k}$ 
$(k=1, 2, ...)$ have the Laurent expansion near $z=0$, 
$$  \frac{\wp_k(\frac{z}{2\pi \sqrt{-1}}, \tau)}{(2\pi \sqrt{-1})^k}
=\frac{1}{z^k}+(-1)^k\sum_{n=1}^{\infty}{2n+1\choose k-1}
E_{2n+2}(\tau)z^{2n+2-k} \eqno{(2.29)}  $$
and the $q$-expansion of $\frac{\wp_1(\frac{z}{2\pi\sqrt{-1}}, \tau)}{2\pi
\sqrt{-1}}$ is
$$ \frac{{\wp}_1(\frac{z}{2\pi \sqrt{-1}}, \tau)}{2\pi \sqrt{-1}}
=E_2(\tau)z+ \hf {e^z+1\over e^z-1}
+ \sum_{n=1}^{\infty}\left({q^n \over e^z-q^n}
-{e^zq^n\over 1-e^zq^n}\right), \eqno{(2.30)}$$
see \cite{La} p.248.

Zhu \cite{Zh} introduced a formal power series. We adopt it by 
multiplying $(2\pi\sqrt{-1})^{-k}$:   
$$ P_k(z, q)={1\over (k-1)!}
\sum_{n\not=0}\frac{n^{k-1}z^n}{1-q^n} \eqno{(2.31)} $$
where ${1\over 1-q^n}$ is understood as $ \sum_{i=0}^{\infty} q^{ni}$ 
and $\frac{1}{1-q^{-n}}=-q^n\frac{1}{1-q^n}$ for $n>0$. 
Note that $z{d\over dz}P_k(z, q)=kP_{k+1}(z, q)$. 
It is easy to prove that $P_k(z, q)$ converges uniformly and 
absolutely in every closed subset of the domain 
$ \{(z, q)|\ |q|<|z|<1\}$. 
The relation between $P_k(z, q)$ and $\wp_k(z, \tau)$ are given by 
$$ P_1(e^z, q)=-\frac{1}{2\pi \sqrt{-1}}\wp_1(\frac{z}{2\pi \sqrt{-1}}, \tau)
+E_2(\tau)z-\frac{1}{2}  \eqno{(2.32)} $$
$$P_2(e^z, q)=-\frac{1}{(2\pi\sqrt{-1})^2}
\wp_2(\frac{z}{2\pi \sqrt{-1}}, \tau)+E_2(\tau)  
\quad \mbox{ and }  \eqno{(2.33)} $$
$$P_k(e^z, q)=(\frac{-1}{2\pi\sqrt{-1}})^k
\wp_k(\frac{z}{2\pi\sqrt{-1}}, \tau)  \mbox{ for }k>2. 
\eqno{(2.34)}$$

\subsection{Equations}
In this paper, we will quote several equations from \cite{Zh}.

Write $(1+z)^{\wt(a)-1}(ln(1+z)^{-1})=\sum_{i\geq -1}c_iz^i$. 
Then $c_{-1}=1$ and 
$$  a[-1]u=\sum_{i\geq -1} c_ia(i)u  \eqno{(2.35)} $$
by the definition of $Y[a, z]$. 
Zhu proved the following equations (c.f. (4.3.8)-(4.3.11) in \cite{Zh}): 
$$ \begin{array}{l}
\sum_{i\geq -1}c_i\Res_w((w-z)^iz^{n-1-i}w^{-n}
-\sum_{i\geq -1}c_i\Res_w((-z+w)^iz^{n-1-i}w^{-n} \cr
=1,  \end{array}  \eqno{(2.36)} $$
$$ \begin{array}{l}
\sum_{i\geq -1}c_i\Res_w((w-z)^iz^{\wt(a)-1-i}w^{-\wt(a)}P_1(\frac{z}{w},q) \cr
-\sum_{i\geq -1}c_i\Res_w((-z+w)^iz^{\wt(a)-1-i}
w^{-\wt(a)}(P_1(\frac{zq}{w},q)-1) \cr
=-\frac{1}{2}  \end{array}  \eqno{(2.37)} $$
and for $m\geq 2$, 
$$ \begin{array}{l}
\sum_{i\geq -1}c_i\Res_w((w-z)^iz^{\wt(a)-1-i}w^{-\wt(a)}P_m(\frac{z}{w},q)\cr
-\sum_{i\geq -1}c_i\Res_w((-z+w)^iz^{\wt(a)-1-i}w^{-\wt(a)}
P_m(\frac{zq}{w}, q) \cr
=E_m(q).  \end{array}  \eqno{(2.38)} $$

He also got the following equation at the end of the proof of Proposition 
4.3.2 in \cite{Zh}. 

$$ \begin{array}{l}
\sum_{i \in \N}\sum_{k=1}^{\infty}
\left( {\wt(a)-1+k\choose i}\frac{1}{1-q^k}
x^k+{\wt(a)-1-k\choose i}\frac{1}{1-q^{-k}}x^{-k} \right)a(i)b \cr
=\sum_{m \in \N}P_{m+1}(x, q)a[m]b.   
\end{array} \eqno{(2.39)}  $$

\pr 
Set ${m-1+x\choose i}=\sum_{s=0}^ic(m,i,s)x^s$.  Then since 
$$\begin{array}{l}
\sum_{s=0}^{\infty}\frac{1}{s!}[ln(1+z)]^s(1+z)^{\wt(a)-1}w^s
=(1+z)^{\wt(a)-1}\exp(ln(1+z)w) \cr
=(1+z)^{\wt(a)-1}(1+z)^w=\sum_{i=0}^{\infty} {\wt(a)-1+w\choose i}z^i \cr
=\sum_{i=0}^{\infty}\sum_{s=0}^i c(\wt(a),i,s)w^sz^i
=\sum_{s=0}^{\infty}\sum_{i=s}^{\infty} c(\wt(a),i,s)z^iw^s, 
\end{array}$$
we have 
$$\sum_{i=s}^{\infty} c(\wt(a),i,s)z^i
=\frac{1}{s!}[ln(1+z)]^s(1+z)^{\wt(a)-1}$$
and so 
$$ \sum_{i=s}^{\infty}c(\wt(a),i,s)a(i)=\frac{1}{s!}a[s]. $$
Therefore 
$$ \begin{array}{l}
\sum_{m=0}^{\infty}\sum_{n\not=0}{\wt(a)-1+n\choose m}
\frac{x^n}{1-q^m}a(m)b=\sum_{m=0}^{\infty}\sum_{n\not=0}
(\sum_{s=0}^m c(\wt(a),m,s)n^s)\frac{x^n}{1-q^n}a(m)b \cr
=\sum_{s=0}^{\infty}\sum_{n\not=0}(\sum_{m=s}^{\infty} c(\wt(a),m,s)a(m)b)
\frac{n^sx^n}{1-q^n} \cr
=\sum_{s=0}^{\infty}\sum_{n\not=0}\frac{1}{s!}a[s]b\frac{n^sx^n}{1-q^n}
=\sum_{s=0}^{\infty}P_{s+1}(x,\tau)a[s]b. 
\end{array} $$

\section{Trace functions}
Let $W$ be a $V$-module. For an intertwining operator 
$I(\ast, z) \in I{W\choose U\quad W}$, 
define \\
$F^I_W:(V\otimes \CH)^{\otimes(i-1)}\otimes (U\otimes \CH)\otimes 
(V\otimes \CH)^{\otimes(n-i)} 
\longrightarrow \C$ by   
$$ \begin{array}{c} 
F^I_W((a_1, z_1), ..., 
(a_{i-1}, z_{i-1}), (u, z_i), (a_{i+1}, z_{i+1}), ..., (a_n, z_n))\cr
=z_1^{{\wt}(a_1)}...z_i^{{\wt}(u)}...z_n^{{\wt}(a_n)}{\tr}|_W
Y^W(a_1, z_1)Y^W(a_2, z_2)...I(u, z_i)...Y^W(a_n, z_n)q^{L(0)}  
\end{array} \eqno{(3.1)} $$
for homogeneous $a_i \in V\ (i=1, ..., n)$ and $u \in U$ and 
extend it linearly. In this paper, we will use $F^I_W$ only 
for $n=1, 2$. 
For simplicity we will often omit the lower index $W$ in $F^I_W$ 
when no confusion should arise. 
For example, $F^I((u, z), \tau)=\tr_{|W} o^I(u)q^{L(0)}$ 
does not depend on $z$ and so we denote it $F^I_W(u, \tau)$.

Since we will calculate the traces and all 
coefficients $a_i(m)$ of $Y(a_i, z)$ shift the grading by integers, 
it is sufficient to consider only coefficients $u(k)$ of 
$I(u, z_i)=\sum u(k)z_i^{-k-1}$ which shift the grading by integers.

The case $U=V$ is Zhu's theory. 
The arguments (circulating arguments) in Zhu's paper depends on 
Commutativity among $ \{Y^W(a_i, z):i\}$, but not on 
Commutativity of $Y^W(a_i, z)$ with itself. 
Hence if only one intertwining operator $I(u, z_i)$ appears 
in the definition of $F_W^I$, 
$ \{Y^W(a_1, z_1), ..., Y^W(a_{i-1}, z_{i-1}), I(u, z_i), ..., 
T^W(a_n, z_n)\}$ 
satisfy Commutativity each other, that is, with the others, 
and so we can apply the circulating arguments. 
Therefore we have the following results. 

\begin{prn}  For any $a \in V$ and $u \in U$, we have 
$$ F^I((a[0]u, z), q)=0. \eqno{(3.2)} $$
\end{prn}

\pr
$$ \begin{array}{l}
0=\tr_{|W}(o(a)I(u, z)q^{L(0)}-I(u, z)q^{L(0)}o(a))=
\tr_{|W}[o(a), I(u, z)]q^{L(0)}\cr
=\tr_{|W} \sum_{i=0}^{\infty}{\wt(a)-1\choose i}z^{\wt(a)-1-i}
I(a(i)u, z)q^{L(0)} \cr
=z^{-\wt(a)}F^I((a[0]u, z), q)
\end{array}$$
\prend

\begin{prn}
$$F^I((a, x), (u, z), q)=z^{-\wt(u)}\tr_{|W}o(a)o^I(u)q^{L(0)}
+ \sum_{m \in \N}P_{m+1}(\frac{z}{x}, q)o^I(a[m]u)q^{L(0)} 
\eqno{(3.3)}$$
$$\begin{array}{l}
F^I((u, z), (a, x), q) \cr
=z^{-\wt(u)}\tr_{|W}o(a)o^I(u)q^{L(0)}
+ \sum_{m \in \N}\left(P_{m+1}(\frac{zq}{x}, q)
-\delta_{m, 0}\right)o^I(a[m]u)q^{L(0)} 
\end{array} \eqno{(3.4)}$$
\end{prn}

\pr
We just follow the proof of Proposition 4.3.2 in \cite{Zh} with 
suitable modifications. 
For $k\not=0$, set $o_k(a)=a(\wt(a)-1+k)$. We have 
$$\begin{array}{l} 
\tr_{|W} o_k(a)I(u, z)q^{L(0)}
=\tr_{|W} [o_k(a), I(u, z)]q^{L(0)}+ \tr_{|W} I(u, z)o_k(a)q^{L(0)} \cr
=\sum_{i \in \N}{\wt(a)-1+k\choose i}z^{\wt(a)-1+k-i}\tr_{|W} I(a(i)u, z)q^{L(0)}
+ \tr_{|W} I(u, z)q^{L(0)}o_k(a)q^k \cr
=\sum_{i \in \N}{\wt(a)-1+k\choose i}z^{\wt(a)-1+k-i}\tr_{|W} 
I(a(i)u, z)q^{L(0)}+ \tr_{|W} o_k(a)I(u, z)q^{L(0)}q^k 
\end{array} $$
Solving for $ \tr_{|W} o_k(a)I(u, z)q^{L(0)}$ in the above identity, 
we have 
$$\tr_{|W}o_k(a)I(u, z)q^{L(0)}
=\frac{1}{1-q^k}\sum_{i \in \N}{\wt(a)-1+k\choose i}
z^{\wt(a)-1+k-i}\tr_{|W}I(a(i)u, z)q^{L(0)}. \eqno{(3.5)}$$
Similarly, we have 
$$\tr_{|W} I(u, z)o_k(a)q^{L(0)}
=\frac{q^k}{1-q^k}\sum_{i \in \N}{\wt(a)-1+k\choose i}
z^{\wt(a)-1+k-i}\tr_{|W}I(a(i)u, z)q^{L(0)}. \eqno{(3.6)}$$

Using the above, we have 
$$ \begin{array}{l}
F^I((a, x), (u, z), q) \cr
=z^{\wt(u)}\tr_{|W} o(a)I(u, z)q^{L(0)}+z^{\wt(u)}\sum_{k\not=0}x^{-k}\tr_{|W} 
o_k(a)I(u, z)q^{L(0)} \cr
=z^{\wt(u)}\tr_{|W} o(a)I(u, z)q^{L(0)} \cr
+z^{\wt(u)}\sum_{k\not=0}x^{-k}\frac{1}{1-q^k}
\sum_{j \in \N}{\wt(a)-1+k\choose i}z^{\wt(a)-1+k-i}
\tr_{|W}I(a(i)u, z)q^{L(0)} \cr
=z^{\wt(u)}\tr_{|W} o(a)I(u, z)q^{L(0)} \cr
+ \sum_{i \in \N}\sum_{k=1}^{\infty}({\wt(a)-1+k\choose i}\frac{1}{1-q^k}
(\frac{z}{x})^k+{\wt(a)-1-k\choose i}\frac{1}{1-q^{-k}}(\frac{z}{x})^{-k})
F^I((a(i)u, z), q) \cr
=\tr_{|W} o(a)o^I(u)q^{L(0)}+ \sum_{m \in \N}
\frac{P_{m+1}(\frac{z}{x}, q)}{(2\pi i)^{m+1}}o^I(a[m]u)q^{L(0)} 
\quad \mbox{ by (2.39)}.
\end{array} $$
Similarly, we have
$$ \begin{array}{l}
F^I((u, z), (a, x), q) \cr
=z^{\wt(u)}\tr_{|W} I(u, z)o(a)q^{L(0)}+z^{\wt(u)}\sum_{k\not=0}x^{-k}\tr_{|W} 
I(u, z)o_k(a)q^{L(0)} \cr
=z^{\wt(u)}\tr_{|W} I(u, z)o(a)q^{L(0)} \cr
+z^{\wt(u)}\sum_{k\not=0}x^{-k}\frac{q^k}{1-q^k}
\sum_{j \in \N}{\wt(a)-1+k\choose i}
z^{\wt(a)-1+k-i}\tr_{|W}I(a(i)u, z)q^{L(0)} \cr
=z^{\wt(u)}\tr_{|W} o(a)I(u, z)q^{L(0)}-F^I((a[0]u, z), q)  \cr
+ \sum_{i \in \N}\sum_{k=1}^{\infty}({\wt(a)-1+k\choose i}\frac{1}{1-q^k}
(\frac{qz}{x})^k+{\wt(a)-1-k\choose i}\frac{1}{1-q^{-k}}(\frac{qz}{x})^{-k})
F^I((a(i)u, z), q) \cr
=\tr_{|W} o(a)o^I(u)q^{L(0)}
+ \sum_{m \in \N}(\frac{P_{m+1}(\frac{qz}{x}, q)}{(2\pi i)^{m+1}}
-\delta_{m, 0})o^I(a[m]u)q^{L(0)} \quad \mbox{ by (2.39)}. 
\end{array} $$

\prend

\begin{prn}
$$\tr_{|W} o(a)o(u)q^{L(0)}
=\tr_{|W} o(a[-1]u)q^{L(0)}-\sum_{k=1}^{\infty}E_{2k}
(\tau)\tr_{|W} o(a[2k-1]u)q^{L(0)} \eqno{(3.7)}$$
\end{prn}

\pr  Write $(1+z)^{\wt(a)-1}(ln(1+z))^{-1}=\sum_{i\geq -1}c_iz^i$ (note that 
$c_{-1}=1$).  Then $a[-1]u=\sum_{i\geq -1} c_ia(i)u$. 
We have 
$$ \begin{array}{l}
F^I((a[-1]u, z), q)
=\sum_i c_iz^{\wt(a)-1-i}z^{\wt(u)}\tr_{|W} I(a(i)u, z)q^{L(0)} \cr
=\sum_i c_iz^{\wt(a)-1-i+ \wt(u)}\Res_{w-z}(w-z)^i
\tr_{|W} I(Y(a, w-z)u, z)q^{L(0)} \cr
=\sum_i c_iz^{\wt(a)-1-i+ \wt(u)}\Res_{w}(w-z)^i
\tr_{|W} Y(a, w)I(u, z)q^{L(0)} \cr
-\sum_i c_iz^{\wt(a)-1-i+ \wt(u)}\Res_{w}(-z+w)^i
\tr_{|W} I(u, z)Y(a, w)q^{L(0)} \cr
=\sum_ii c_i\Res_w( (w-z)^iz^{\wt(a)-1-i}w^{-\wt(a)}
F^I((a, w), (u, z), q)) \cr
-\sum_ii c_i\Res_w( (-z+w)^iz^{\wt(a)-1-i}w^{-\wt(a)}
F^I((u, z), (a, w), q)).
\end{array} $$
By Proposition 3.2 and (2.36)-(2.38), we have 
$$\begin{array}{l}
o^I(a[-1]u)q^{L(0)}
=z^{\wt(u)}\tr_{|W} o(a)I(u, z)q^{L(0)} \cr
-\hf F^I((a[0]u, z), q)+ \sum_{k=1}^{\infty}E_{2k}(\tau)
F((a[2k-1]u, z), q) \cr
=\tr_{|W} o(a)o^I(u)q^{L(0)}
+ \sum_{k=1}^{\infty}E_{2k}(\tau)o^I(a[2k-1]u)q^{L(0)}. 
\end{array}$$
\prend 

As an intertwining operator version of Proposition 4.3.6 in \cite{Zh}, 
we have proved the following: 

\begin{thm} For $a \in V, u \in U$, we have 
$$ {\tr}|_Wo(a[0]u)q^{L(0)}=0,     \eqno{(3.8)} $$
$$ {\tr}|_Wo(a[-2]u)q^{L(0)}+ \sum_{k=2}^{\infty}
(2k-1)E_{2k}(\tau){\tr}|_Wo(a[2k-2]u)q^{L(0)}=0.  \eqno{(3.9)}$$ 
\end{thm}

\pr
(3.8) is proved in Proposition 3.1. 
Replace $a$ in (3.7) by $L[-1]a$. Since \\
$(L[-1]a)[2k-1]=-(2k-1)a[2k-2]$ 
and $o(L[-1]a)=0$, we have $(3.9)$. 
\prend

\begin{lmm}  For every $u \in U$, 
$$\tr_{|W} o(\tilde{\omega})o(u)q^{L(0)}
=(q\frac{d}{dq}-\frac{c}{24})
\tr_{|W} o(u)q^{L(0)} \eqno{(3.10)}$$
and so 
$$\tr_{|W} o(L[-2]u)q^{L(0)}
-\sum_{k=1}^{\infty}E_{2k}(\tau)\tr_{|W} o(L[2k-2]u)
q^{L(0)}=(q\frac{d}{dq}-c/24)\tr_{|W} o(u)q^{L(0)} \eqno{(3.11)}$$
\end{lmm}

\pr Clearly, $ \tr_{|W} o(\omega)o(u)q^{L(0)}=\tr_{|W} L(0)o(u)q^{L(0)}
=(q\frac{d}{dq})\tr_{|W} o(u)q^{L(0)}$. 
Since \\
$o(\tilde{\omega})=o(\omega)-\frac{c}{24}$, we have (3.10).  
Substitute $ \tilde{\omega}$ into $a$ of (3.7), we have (3.11) by (3.10). 
\prend

\section{The space of one point functions on the torus}
We will use the following notation: \\
(a)  $M(\Gamma(1))$ is the ring of holomorphic modular forms on 
$ \Gamma(1)=SL_2(\Z)$; it is naturally graded 
$M(\Gamma(1))=\oplus_{k\geq 0}M_k(\Gamma(1))$, where 
$M_k(\Gamma(1))$ is the space of forms of weight $k$. 
It is known that 
$M(\Gamma(1))$ is generated by $E_4(\tau)$ and $E_6(\tau)$. 
(c.f. Proposition 1.3.4 in \cite{Bu}.)\\
(b) Set $U(\Gamma(1))=M(\Gamma(1))\otimes_{\C} U$. \\
(c) $O_q(U)$ is the $M(\Gamma(1))$-subspace of $U(\Gamma(1))$ generated by 
the following elements:
$$  v[0]u, \quad v \in V, u \in U   \eqno{(4.1)} $$
$$  v[-2]u+ \sum_{k=2}^{\infty} (2k-1)E_{2k}(\tau)\otimes v[2k-2]u  
  \quad v \in V, u \in U \eqno{(4.2)} $$

A crucial connection between Zhu-algebra and Eisenstein series is 
Lemma 5.3.2 in \cite{Zh}.  We will reform it for modules as follows: 

\begin{lmm} 
Set $a*_{\tau}u=a[-1]u-\sum_{k=1}^{\infty}E_{2k}(\tau)a[2k-1]u$.  
The constant terms of $a*_{\tau}u$ in $U[[q]]$ for $a \in V, u \in U$ is 
$$a\cdot u-\hf a[0]u. \eqno{(4.3)} $$ 
In particular, the constant term of 
$ \tilde{\omega}*_{\tau}u$ $(=L[-2]u-\sum_{i=1}^{\infty}E_{2k}(\tau)
L[2k-2]u)$ is 
$$\tilde{\omega}\cdot u-\hf \tilde{\omega}[0]u. \eqno{(4.4)} $$ 
\end{lmm}

\pr
By (2.29), we have 
$a[-1]b-\sum_{k=1}^{\infty}E_{2k}(\tau)a[2k-1]u
=\Res_z(Y[a, z](\wp_1(\frac{z}{2\pi i}, \tau)
-\frac{G_2(\tau)}{(2\pi i)^2}z))u$ and its constant term is 
$$ \begin{array}{l}
\hf\Res_z\left( Y[a, z]\frac{e^z+1}{e^z-1}\right)u \cr
=\hf\Res_z\left( Y(a, e^z-1)e^{z \wt(a)}
\frac{e^z+1}{e^z-1}\right)u \cr
=\hf\Res_w\left( Y(a, w)\frac{(1+w)^{\wt(a)}}{w}
\frac{(w+2)}{(1+w)}\right)u \cr
 =a\cdot u-\hf a[0]u \qquad \mbox{ by (2.15) and (2.4)} \end{array} $$  
\prend

\noindent
\begin{dfn}
For an irreducible $V$-module $U$, 
we now define the space $C_1(U)$ of one point functions on $U$ 
to be the $ \C$-linear space consisting of functions
$$  S: U(\Gamma(1))\otimes \CH \to \C   $$
satisfying the following conditions: \\
(C1)  For $u \in U(\Gamma(1))$, $S(u, \tau)$ is holomorphic in $ \tau$. \\
(C2)  $S(u, \tau)$ is $M(\Gamma(1))$-linear 
in the sense that $S(u, \tau)$ is $ \C$-linear 
in $u$ and satisfies \\
$$ S(f(\tau)\otimes u, \tau)=f(\tau)S(u, \tau)   $$
for $f(\tau) \in M(\Gamma(1))$ and $u \in U$. \\
(C3) For $u \in O_q(U)$, $S(u, \tau)=0$  \\
(C4) For $u \in U$, 
$$  S(L[-2]u, \tau)={\pd} S(u, \tau)+
\sum_{k=2}^{\infty}E_{2k}(\tau)S(L[2k-2]u, \tau)$$
Here $ \pd$ is the operator which is linear in $u$ and satisfies 
$$ \pd S(u, \tau)=\pd_k S(u, \tau)=\frac{1}{2\pi i}\frac{d}{d\tau}S(u, \tau)
+kE_2(\tau)S(u, \tau) \eqno{(4.5)}$$
for $u \in U_{[k]}$ and 
$$ \pd S(f\otimes u, \tau)=(\pd_hf(\tau))S(u, \tau)+f(\tau)\pd S(u, \tau)  $$
for $f(\tau) \in M_h(\Gamma(1))$ and 
$ \pd_hf(\tau)=\frac{1}{2\pi i}\frac{df(\tau)}{d\tau}+hE_2(\tau)f(\tau)$. 
We note $q\frac{d}{dq}=\frac{1}{2\pi i}\frac{d}{d\tau}$. 
\end{dfn}

Set $S^I(u, \tau)=F^I(u, \tau)q^{-c/24}=\tr_{|W}o^I(u)q^{L(0)-c/24}$ 
for $u \in U$ and extend it 
linearly for $M(\Gamma(1))\otimes U$, where $c$ is the central 
charge of $V$.

\begin{prn}
  $$S^I(\ast, \tau) \in C_1(U) \mbox{  if $S^I(u, \tau)$ is 
holomorphic in $ \tau$}.$$
\end{prn}

\pr
(C2) is clear. 
Theorem 3.4 implies (C3). By Lemma 3.5, 
$$\begin{array}{l}
 F^I(L[-2]u, \tau)q^{-c/24}
-\sum_{k=2}^{\infty}E_{2k}(\tau)F^I(L[2k-2]u, \tau)q^{-c/24} \cr
=E_2(\tau)F^I(L[0]u, \tau)q^{-c/24}+(q\frac{d}{dq}-c/24)
\tr_{|W}o(u)q^{L(0)})q^{-c/24} \cr
=E_2(\tau)F^I(L[0]u, \tau)q^{-c/24}+q\frac{d}{dq}
\tr_{|W}o(u)q^{L(0)-c/24} \cr
=\pd \left( \tr_{|W} o(u)q^{L(0)-c/24}\right) \mbox{ by (4.5)}. 
\end{array}$$
Hence $S^I(\ast, \tau)$ satisfies (C4). 
\prend

We fix an element $S(\ast, \tau) \in C_1(U)$ for a while. 

\begin{lmm}
Let $u \in U$. If $U$ satisfies condition $C_{[2,0]}$ 
then there are $m \in \N$ and 
$r_i(\tau) \in M(\Gamma)$ for $i=0, ..., m-1$ such that 
$$ S(L[-2]^mu, \tau)+ \sum_{i=0}^{m-1}r_i(\tau)S(L[-2]^iu, \tau)=0. 
\eqno{(4.6)}$$
\end{lmm}

\pr  
Since $C_{[2,0]}(U)$ is a direct 
sum of homogeneous subspaces with respect to the degree $[\wt]$, there exists 
$N$ such that elements $a$ satisfying $[\wt]a>N$ are in $C_{[2,0]}(U)$.  
Let $A$ be the 
$M(\Gamma(1))$-submodule of $U(\Gamma(1))$ generated by 
$ \oplus_{n\leq N}U[n]$. 
We claim that $U=A+O_q(U)$. 
Suppose false and let $K$ be a minimal weight of $U/(A+O_q(U))$ and 
choose $u \in U-(A+O_q(U))$ with weight $K=[\wt]u$.
Clearly, we have $[\wt]u>N$ and so $u=\sum b_i[-2]w_i+\sum c_j[0]u_j$ with 
$b_i, c_j \in V$, $w_i,u_j \in U$. We may assume that $u=b[-2]w$ and 
$[\wt](b[-2]u)=K$. Then 
$$b[-2]w \in -\sum_{k=2}^{\infty}(2k-1)E_{2k}(\tau)b[2k-2]w+O_q(U) 
\subseteq A+O_q(U)   \eqno{(4.7)}$$ 
by the minimality of $K$ since $[\wt](b[2k-2]w)< K$ for 
$k\geq 1$. Hence $U(\Gamma(1))/Q_q(V)$ is a finitely generated 
$M(\Gamma(1))$-module and so 
there are $m \in \N$ and $r_i(\tau) \in M(\Gamma)$ such that 
$L[-2]^mu+ \sum_{i=0}^{m-1}r_i(\tau)L[-2]^iu \in O_q(U)$. 
\prend

\begin{lmm} Let $u \in U$. 
Suppose that $L[n]u=0$ for $n>0$.  Then for $i\geq 1$ there are 
elements $f_j(\tau) \in M(\Gamma)$ for $j=0, ..., i-1$ such that 
$$ S(L[-2]^iu, \tau)=\pd^iS(u, \tau)+ \sum_{j=0}^{i-1}f_j(\tau)\pd^jS(u, \tau) 
\eqno{(4.8)}$$
\end{lmm}

\pr
We will prove (4.8) by induction on $i$. Since $L[n]u=0$ for $n>0$, 
(C4) implies $S(L[-2]u, \tau)=\pd S(u, \tau)$. By (C4), we have 
$S(L[-2]L[-2]^iu, \tau)=\pd S(L[-2]^iu, \tau)
+ \sum_{k=2}^{\infty}E_{2k}S(L[2k-2]L[-2]^iu, \tau)$. On the other hand, since 
$L[2k-2]L[-2]^iu$ is written as 
$ \sum_{j=0}^{i-k-1} p_jL[-2]^ju$ with some constants $p_j$ for 
$k\geq 2$, we have the desired statement. 
\prend

Bearing in mind the definition of $ \pd$ and using (4.8), 
(4.6) may be reformulated as follows:

\begin{prn}
Suppose that $U$ satisfies condition $C_{[2,0]}$ and $u \in U$ satisfies 
$L[n]u=0$ for $n>0$.  Then there are $m \in \N$ and $h_i(\tau) \in M(\Gamma)$ 
for $i=1, ..., m-1$ such that 
$$  (\frac{d}{d\tau})^mS(u, \tau)
+ \sum_{i=0}^{m-1}h_i(\tau)(\frac{d}{d\tau})^iS(u, \tau)
=0.  \eqno{(4.9)} $$
\end{prn}

\begin{thm}  If $U$ satisfies condition $C_{[2,0]}$, then 
$S^I(u,\tau)$ is holomorphic on the upper half-plane for $u\in U$. 
\end{thm}

\pr By (3.8) and (3.9), we may assume that $L[n]u=0$ for $n>0$. 
Then by the same argumets as in the proof of Proposition 4.6, we obtain 
that $S^I(u,\tau)$ satisfies (4.9).  Since (4.9) is regular 
and $h_i(\tau)$ converges absolutely and uniformly on every 
closed subset of $\{q\ |\ |q|<1\}$, so does $S^I(u,\tau)$. 
\prend

As a corollary, we have: 

\begin{cry}  Assume that $U$ satisfies condition $C_{[2,0]}$.  
Then $S(u, \tau)$ 
converges absolutely and uniformly in every closed subset of the domain 
$ \{q\ |\ |q|<1\}$ for every $u \in U$ and the limit function can be written 
as $q^hf(q)$, where $f(q)$ is some analytic function in $ \{q\ |\ |q|<1\}$. 
\end{cry}

\pr
We will prove the assertion by induction on $[\wt]u$. 
Assume $L[n]u=0$ for $n>0$.  
Since $h_i(q)$ in (4.9) converges absolutely and uniformly on every closed 
subset of $ \{q\ |\ |q|<1\}$ and the equation (4.9) is regular, we 
have the desired result for $S(u, \tau)$. If $u=L[-1]w$, then 
$S(u, \tau)=0$. So we may assume that $u=L[-2]w$ since $<L[n] : (n<0)>$ is 
generated by $L[-1]$ and $L[-2]$. 
In this case, since 
$$S(L[-2]w, \tau)=\pd S(w, \tau)+ \sum_{k=1}^{\infty}E_{2k}S(L[2k-2]w, \tau), 
$$
we have the desired result by induction. 
\prend

Now assume that $u \in U$ satisfies $L[n]u=0$ for $n>0$. 
(4.9) is a homogeneous linear differential equation with holomorphic 
coefficients and, such that $0$ is a regular singular point. We are 
therefore in a position to apply the theory of Frobenius-Fuchs 
concerning the nature of the solutions to such equations. 
Frobenius-Fuchs theory tells us that $S(u, \tau)$ may be expressed in the 
following form: for some $p\geq 0$, 
$$  S(u, \tau)=\sum_{i=0}^p(\log q)^iS_i(u, \tau)   \eqno{(4.10)}$$
where 
$$  S_i(u, \tau)=\sum_{j=1}^{b(i)}q^{\la_{ij}}S_{i, j}(u, \tau) 
\eqno{(4.11)} $$
$$  S_{i, j}(u, \tau)=\sum_{n=0}^{\infty} a_{i, j, n}q^n  \eqno{(4.12)}$$
is homomorphic on the upper half-plane, and 
$ \la_{i, j_1}\not\equiv \la_{i, j_2} \pmod{\Z}$ for $j_1\not=j_2$.  

We claim that $S(u, \tau)$ has the similar form for every $u \in U$. 
(c.f. Theorem 6.5 \cite{DLiM}).

\begin{thm}  Suppose that $U$ satisfies condition $C_{[2,0]}$. 
For every $u \in U$, 
the function $S(u, \tau)$ can be expressed in the form (4.10), (4.11), (4.12). 
\end{thm}

\pr
We will prove the theorem by induction on $[\wt]u$. If $L(n)u=0$ 
for $n>0$, then we have already mentioned at (4.10)-(4.12). 
Since $S(L[-1]u, \tau)=0$, we may assume $u=L[-2]w$.  Then 
$$S(L[-2]w, \tau)=\pd S(w, \tau)
+ \sum_{k=1}^{ \infty}E_{2k}(\tau)S(L[2k-2]w, \tau). \eqno{(4.13)}$$ 
By induction, $ \pd S(w, \tau)$ and $S(L[2k-2]w, \tau)$ have the desired 
forms and so does \\
$S(L[-2]w, \tau)$. 
\prend

Before we will prove the main assertion, we will have the following lemmas: 

\begin{lmm} Let $A$ be a semi-simple associative algebra over ${\C}$, 
let $ \omega$ be in the center of $A$, and  
let ${}_AB_A$ be an $A$-bimodule and 
let $F$ be a linear functional of $B$ satisfying 
$F(ab)=F(ba)$ for every $a \in A, b \in B$.  
Assume further that $F((\omega-r)^Nu)=0$ for every $u \in U$, where 
$r$ is constant and $N$ is an integer.
Then there are irreducible $A$-modules $M^1, ..., M^n$ on which 
$ \omega$ acts as a scalar $r$ 
and $I^i(\ast) \in \Hom_A(B\otimes_A M^i, M^i)$ for $i=1, ..., n$ such that 
$$  F(b)=\sum_{i=1}^n {\tr}|_{M^i}P_{I^i}(b)  \eqno{(4.14)} $$
for every $b \in B$, where 
$P_{I^i}(b) \in End(M^i)$ is given by $P_{I^i}(b)m_i=I^i(b\otimes m_i)$ 
for $m_i \in M^i$. 
\end{lmm}

\pr 
Since $A$ is a semi-simple associative algebra, $A$ decomposes into 
the direct sum $ \oplus A_i$ of matrix rings $A_i$. 
We may assume that $B$ is an irreducible $A$-bimodule. 
If $A_iB=0$ or $BA_i=0$ for all $i$ then we have 
$F(ab)=F(ba)=F(0)=0$ for all $a \in A_i$ and so  
$F(b)=0 \forall b \in B$. In particular the assertion is held. 
So we may assume that $A$ is a matrix ring and there 
is a bimodule isomorphism $f:{}_AA_A\cong {}_AB_A$ given by 
$f(a)=af(1)$. 
Set $G(a)=F(af(1))$. Since 
$G(aa')=F(aa'f(1))=F(af(a'))=F(f(a'1)a)=F(a'f(1)a)=F(a'f(a))=
F(a'af(1))=G(a'a)$, there is $k \in \C$ such that $G(a)=k\tr(a)$. 
Let ${}_AM$ be an irreducible $A$-module and fix $0\not=d \in M$.
Set $I(af(1)\otimes a'd)=kaa'd$. Then 
$I(\ast) \in \Hom_A(B\otimes_AM, M)$ and 
$P_I(af(1))a'm=I(af(1)\otimes a'm)=kaa'm$ and so 
$ \tr_{|M} P_I(af(1))=k\tr_{|M}a=G(a)=F(af(1))$. Since $ \omega$ is in the 
center of $A$, $ \omega$ acts on $M$ as a scalar, which should be $r$. 
\prend

\begin{lmm}
We have: 
$$ S_{k, j}(a[0]u, \tau)=0 \mbox{ for any } a \in V, u \in U, k, j. 
\eqno{(4.15)}$$
$$ S_{k, j}(u_q, \tau)=0 \mbox{ for any } u_q \in O_q(U), k, j.   
\eqno{(4.16)}  $$
$$ S_{p, j}(\tilde{\omega}*_{\tau}u, \tau)
=\pd S_{p, j}(u, \tau) \mbox{ for }u \in U.  \eqno{(4.17)} $$
$$ S_{k-1}(\tilde{\omega}*_{\tau}u, \tau)
=kS_k(u_q, \tau)+ \pd S_{k-1}(u, \tau) \mbox{ for } u \in U.  
\eqno{(4.18)} $$
$$  (\tilde{\omega}*_{\tau}-\pd)^NS_{p-k, j}(u, \tau)=0 
\mbox{ for }N\!>\!>\!0.   \eqno{(4.19)} $$
\end{lmm}

\pr
Since $U(\Gamma(1))\subseteq U[[q]]$, (4.15) and (4.16) are clear. 
Since $ \tilde{\omega}*_{\tau}u
=L[-2]u-\sum_{k=1}^{\infty}E_{2k}(\tau)L[2k-2]u$, 
we have 
$$\begin{array}{l}
S(\tilde{\omega}*_{\tau}u, \tau)=\pd S(u, \tau) \cr
=\sum_{k=0}^p\{ k(\log q)^{k-1}S_k(u, \tau)+(\log q)^k\pd S_k(u, \tau)\}
\end{array}\eqno{(4.20)} $$
for $u \in U$ and so we have (4.18). 
In particular, we have (4.17). By (4.18), we have 
$(\tilde{\omega}*_{\tau}-\pd)S_{r}(u, \tau)=(r+1)S_{r+1}(u, \tau)$.  
Repeating these steps, we obtain (4.19). 
\prend

\begin{lmm}  Assume $T(u, \tau)=q^{\la}\sum_{n=0}^{\infty}\al_n(u)q^n$ 
satisfies the conditions (4.15), (4.16), \\
(4.19), then the coefficient $ \al_0(u)$ of leading term satisfies 
$$ \al_{0}(a\cdot u)=\al_{0}(u*a) \mbox{ for $u \in U$, $a \in V$, } 
\eqno{(4.21)}$$
$$\al_{0}(u)=0  \mbox{ for $u \in O(U)$ and} \eqno{(4.22)} $$
$$\al_{0}((\omega-c/24-\la)^N\cdot u)=0 \mbox{ for $u \in U$.} \eqno{(4.23)}$$
\end{lmm}

\pr 
Since $T(a[0]u, \tau)=0$ and $a\cdot u-u*a=a[0]u$, (4.21) holds. 
Since the constant terms of elements in $O_q(U)$ 
generates $O(U)$, we have 
$ \al_{0}(u)=0$  for $u \in O(U)$, which proves (4.22). 
Since the constant term of $ \tilde{\omega}*_{\tau}u$ 
is $ \tilde{\omega}\cdot u-\hf\tilde{\omega}[0]u$ for $u \in U$ by 
Lemma 4.1 and $T(\tilde{\omega}[0]u, \tau)=0$ by Proposition 3.1, 
the leading terms of $(\tilde{\omega}*_{\tau}-\pd)T(u, \tau)$ 
is $q^{\la}\al_0((\omega-c/24-\la)\cdot u)$. Since the operator 
$ \tilde{\omega}*_{\tau}-\pd$ does not decrease the minimal degree of  
elements in $U(\Gamma(1))[[q]]$, the leading term of  
$$  (\tilde{\omega}*_{\tau}-\pd)^NT(u, \tau)=0  $$
is $q^{\la}\al_{0}((\omega-c/24-\la)^N\cdot u)=0$ for $u \in U$.  
\prend

We note that $S^{I}(u, \tau)$ satisfies the same conditions.

\begin{lmm}  Assume $T(u, \tau)=q^{\la}\sum_{n=0}^{\infty}\al_n(u)q^n$ 
satisfies the conditions (4.15), (4.16) and (4.19), then there are 
irreducible $V$-modules $ \{W^j:j \}$ with minimal weight $ \la+c/24$ 
and intertwining operators $ \{I^j:j\}$ of type ${W^j \choose U\quad W^j}$ 
such that $T(u, \tau)$ is a sum of $S^{I^j}(u, \tau)$. 
\end{lmm}

\pr
Since $ \al_0(u)$ satisfies (4.21), (4.22), (4.23), Lemma 4.10 implies that 
there are irreducible $A(V)$-modules $W^1(0), ..., W^m(0)$ on which $ \omega$ acts as a 
scalar $c/24+ \la$ and 
$I^i_0 \in \Hom_{A(V)}(A(U)\otimes_{A(V)} W^i(0), W^i(0))$ 
such that 
$ \al_0(u)=\sum_{i=1}^m \tr_{|M^i(0)}P_{I^i_0}(u)$, where 
$P_{I^i_0}(u)m_i=I^i_0(u\otimes m_i)$. 
By Theorem 2.10, there are irreducible $V$-modules $M^i$ with minimal weight 
$ \la$ and 
intertwining operators $I^i(\ast, z) \in I\pmatrix{M^i\cr U\quad M^i}$ 
such that $M^i(0)$ is the top levels of $M^i$ and 
$I^i_0(u\times m_i)=o^{I^i}(u)m_i$ for $m_i \in M^i(0)$. 
Then $T(u, \tau)-\sum S^{I^i}(u, \tau)$ satisfies the same conditions 
(4.15), (4.16) and (4.19) and the degree of the leading term is 
greater than $ \la$. 
Repeating this steps, we finally have the desired result, since 
there are only finitely many non-isomorphic $V$-modules. 
\prend 

In particular, we have 

\begin{prn}  $S_i(\ast, \tau)=0$ if $i>0$. 
\end{prn}

\pr
Suppose $p\geq 1$. By Lemma 4.13, 
$S_{p-1}(\ast, \tau)$ is a linear combination of $S^{I^j}(\ast, \tau)$ and 
so $ \pd S_{p-1}(u, \tau)=S_{p-1}(\tilde{\omega}*_{\tau}u, \tau)$. 
On the other hand, we have 
$S_{p-1}(\tilde{\omega}*_{\tau}u, \tau)
=pS_p(u, \tau)+ \pd S_{p-1}(u, \tau)$ by (4.18), so that 
$S_p(u, \tau)=0$ for $u \in U$.
\prend

So we have proved the following main theorem, which is an intertwining 
operator version of Theorem 5.3.1 in \cite{Zh}. 

\begin{thm}  Suppose $U$ satisfies condition $C_{[2,0]}$. Let 
$ \{W^1, ..., W^m\}$ be the set of irreducible $V$-modules and 
let $ \{ I^{kj}(\ast, z):j=1, ..., j_k \}$ be a basis 
of $I\pmatrix{W^k \cr U\qquad W^k}$.
Then $C_1(U)$ is spanned by 
$$ \{   S_{I^{kj}}(u, \tau) \qquad :k=1, ..., m, j=1, ..., j_k \}.  $$
\end{thm}

\section{Modular invariance}
In this section, we will show that $C_1(U)$ is invariant under the 
action of $SL(2, \Z)$ if the weights of $U$ are real. 
We note that $SL(2, \Z)$ is generated by $ \pmatrix{0&1 \cr -1&0}$ and 
$ \pmatrix{1&1\cr 0&1}$. 
Since $C_1(U)$ is spanned by $S^{I^j}(\ast, \tau)$, which have forms 
$q^r(\sum_{n=1}^{\infty}a_nq^n)$, 
the transformation of $S(u, \tau)$ by 
$ \pmatrix{1&1\cr 0&1}$ is clear and so it is sufficient 
to prove the assertion for $ \ga\tau=\frac{-1}{\tau}$. 

Let $t+2m$ $(-1<\leq t\leq 1, m\in \Z)$ be the lowest weight of $U$. 
Since $U$ is an irreducible $V$-module, 
the weights of elements in $U$ are all in $\Z+t$. 
First we take the principal branch of $(-\imath\tau)^{-t}$ on 
$ \imath \CH$ by taking 
$(re^{2\pi \theta})^{-t}=r^{-t}e^{2\pi -t\theta}$ for $r\geq 0$ and 
$-\frac{1}{2}< \theta\leq \frac{1}{2}$ and $\tau^{-t-2n}$ is understood to be 
$(\tau^{-t})\tau^{-2n}$ for $n \in \Z$.

\begin{thm}
For $S(\ast, \tau) \in C_1(U)$ and 
$ \ga=\pmatrix{0&-1\cr 1&0} \in SL(2, \Z)$ 
define 
$$  S|\gamma(u, \tau)
={(-\imath\tau)}^{-t}\tau^{t-k}S(u, \frac{-1}{\tau})  \eqno{(5.2)}  $$
for $u \in U_{[k]}$.  Then $S|\gamma(\ast, \tau) \in C_1(U)$.
\end{thm}

\pr
Clearly, 
$S|\gamma(u, \tau)
=(-\imath\tau)^{-t}\tau^{t-k}S(u, \frac{-1}{\tau})$ is 
clearly holomorphic and 
it is also clear that $S|\gamma(\ast, \tau)$ satisfies (C2). 
Let $a \in V_{[p]}$. Then we have: 
$$S|\gamma(a[0]u, \tau)
=(-\imath\tau)^{-t}\tau^{t-k-p-1}S(a[0]u, \frac{-1}{\tau})
=0$$
and  
$$ \begin{array}{l}
S|\gamma(a[-2]u
+\sum_{j=2}^{\infty}(2j-1)E_{2j}(\tau)\otimes a[2j-2]u, \tau)\cr
=(-\imath\tau)^{-t}\tau^{t-k-p-1}S(a[-2]u, \tau) \cr
\mbox{}\qquad +(-\imath\tau)^{-t}\sum_{j=2}^{\infty}\tau^{t-k-p+2j-1}
S((2j-1)E_{2j}(\frac{-1}{\tau})\otimes a[2j-2]u, 
\frac{-1}{\tau}) \cr
=(-\imath\tau)^{-t}\tau^{t-k-p-1}S(a[-2]u, \tau) \cr
\mbox{}\qquad 
+(-\imath\tau)^{-t}\tau^{t-k-p-1}\sum_{j=2}^{\infty}S((2j-1)E_{2j}(\tau)\otimes a[2j-1]u, 
\frac{-1}{\tau}) \cr
=(-\imath\tau)^{-t}\tau^{t-k-p-1}S\left(a[-2]u+ \sum_{j=2}^{\infty}(2j-1)
E_{2j}(\tau)\otimes a[2j-1]u, \frac{-1}{\tau}\right) \cr
=0. 
\end{array} $$
We hence have (C3). We also have 
$$ 
\begin{array}{l}
\frac{d}{d\tau}(S|\gamma(u, \tau))
=\frac{d}{d\tau}\left((-\imath\tau)^{-t}\tau^{t-k}S(u, \frac{-1}{\tau})\right)
\cr
=-k(-\imath\tau)^{-t}\tau^{t-k-1}S(u, \frac{-1}{\tau}) 
+(-\imath\tau)^{-t} \tau^{t-k}\frac{d}{d\tau}S(u, \frac{-1}{\tau}). 
\end{array}  $$
Using this, we obtain: 
$$\begin{array}{l}
S|\gamma(L[-2]u, \tau)\cr
=(-\imath\tau)^{-t}\tau^{t-k-2}S(L[-2]u, \frac{-1}{\tau}) \cr
=(-\imath\tau)^{-t}\tau^{t-k-2}\left(
\frac{d}{2\pi id(-1/\tau)}S(u, \frac{-1}{\tau})
+kE_2(\frac{-1}{\tau})S(u, \frac{-1}{\tau}) \right. \cr
\mbox{}\quad + 
\left. \sum_{n=2}^{\infty}E_{2k}(\frac{-1}{\tau})
S(L[2n-2]v, \frac{-1^{}}{\tau_{}}) \right)\cr
=(-\imath\tau)^{-t}\tau^{t-k-2}\left(
\frac{\tau^2}{2\pi i}\frac{d}{d\tau}S(u, \frac{-1}{\tau})
+k\tau^2E_2(\tau)S(u, \frac{-1}{\tau})
-\frac{k\tau}{2\pi i}S(u, \frac{-1}{\tau})
\right. \cr
\mbox{}\quad \left. + \tau^{2n}E_{2n}(\tau)S(L[2n-2]u, \frac{-1^{}}{\tau_{}}) 
\right) \cr
=(-\imath\tau)^{-t}\left\{
\tau^{t-k}\frac{1}{2\pi i}\frac{d}{d\tau}S(u, \frac{-1}{\tau})
+ \tau^{t-k}kE_2(\tau)S(u, \frac{-1}{\tau})
-\tau^{t-k-2}\frac{k\tau}{2\pi i}S(u, \frac{-1}{\tau}) \right. \cr
\mbox{}\left. \quad + \tau^{2n+t-k-2}E_{2n}(\tau)S(L[2n-2]u, \frac{-1}{\tau})
\right\} \cr
=\frac{1}{2\pi i}\left\{ 
(-\imath\tau)^{-t}\tau^{t-k}\frac{d}{d\tau}S(u, \frac{-1}{\tau})
-(-\imath\tau)^{-t}\tau^{t-k-2}k\tau S(u, \frac{-1}{\tau})\right\} \cr 
\mbox{}\quad +(-\imath\tau)^{-t}\tau^{t-k}kE_2(\tau)S(u, \frac{-1}{\tau}) 
+ (-\imath\tau)^{-t}\tau^{2n+t-k-2}E_{2n}(\tau)S(L[2n-2]u, \frac{-1}{\tau}) \cr
=\frac{d}{2\pi id\tau}S|\gamma(u, \tau)
+kE_2(\tau)S|\gamma(u, \frac{-1}{\tau})
+E_{2n}(\tau)S|\gamma(L[2n-2]u, \frac{-1}{\tau}) 
\end{array}$$
and so (C4). 
This completes the proof of Theorem 5.1. 
\prend

\begin{lmm} 
For $v\in U_{[k]}$,  
$$ P:\left\{
\begin{array}{lll}
\pmatrix{0&1 \cr -1&0} & \rightarrow & 
(-\imath\tau)^{-t}\tau^{t-k}S(u,\frac{-1}{\tau}) \cr
\pmatrix{1&1\cr 0&1} &\rightarrow & (\imath)^{t/3}S(u,\tau+1)  
\end{array} \right.  $$
is a representation of $SL(2,\Z)$. 
\end{lmm}

\pr 
Set $A=\pmatrix{0&1\cr -1&0}$ and $B=\pmatrix{1&1\cr 0&1}$. 
It is sufficient to prove $P(A)^2=1$ and $(P(A)P(B))^3=1$ since 
$SL(2,\Z)$ is generated by two elements satisfying these relations freely. 
$$ \begin{array}{l}
(P(A))^2S(u,\tau) \cr
=(-\imath\tau)^{-t}\tau^{t-k}(-\imath\frac{-1}{\tau})^{-t}
(\frac{-1}{\tau})^{t-k}S(u,\frac{-1}{\frac{-1}{\tau}})\cr
=(-\imath\tau)^{-t}(-\imath\frac{-1}{\tau})^{-t}(-1)^{t-k}S(u,\tau) \cr
=S(u,\tau) \qquad \mbox{(since $k-t\in 2\Z$).} 
\end{array}  $$
$$\begin{array}{l}
(P(A)P(B))^3S(u,\tau) \cr
=(\imath)^{t} (-\imath\tau)^{-t} \tau^{t-k}
(-\imath (\frac{-1}{\tau}+1) )^{-t}( \frac{-1}{\tau}+1 )^{t-k}
(-\imath( \frac{-1}{\frac{-1}{\tau} +1}+1))^{-t}
(\frac{-1}{ \frac{-1}{\tau}+1}+1)^{k-m} \times \cr
\mbox{}\qquad \times S(u,\frac{-1}{ \frac{-1}{\frac{-1}{\tau}+1}+1}+1) 
\cr
=(\imath)^t(1-\tau)^{-t}(-\imath\frac{1}{\tau-1})^{-t}S(u,\tau) \cr
=S(u,\tau)
\end{array} $$

\section{Examples}
Let $L(c,h)$ denote the irreducible module of Virasoro algebra 
with the highest weight $h$ and the central charge $c$ and it was 
proved in \cite{FZ} that $L(c,0)$ is a VOA.
The work in \cite{FQS} and \cite{GKO} gives a complete classification 
of unitary highest weight representations of the Virasoro algebra. 
In particular, $L(c_m,0)$ for 
$c=c_m=1-\frac{6}{(m+2)(m+3)} \quad (m=0,1,2,...)$ 
are rational VOAs called discrete series and their 
irreducible modules are $L(c_m,h^m_{r,s})$ with  
$h=h^m_{r,s}=\frac{[(m+3)r-(m+2)s]^2-1}{4(m+2)(m+3)} \quad 
(r,s\in {\N}, 1\leq s\leq r\leq m+1)$.
The fusion rules among $L(c_m,h^m_{r,s})$ are all determined, 
see \cite{FF}, \cite{W}. 
 
\begin{lmm}  $L(c_m,h^m_{r,s})$ satisfies condition $C_{[2,0]}$. 
\end{lmm}

\pr
Set $U=L(c_m,h^m_{r,s})$ and 
let $e$ be a highest weight vector of $U$.  
We note that 
$(L(c,0), Y[,],\1,\omega)\cong (L(c,0), Y(,),\1,\omega)$ 
as VOAs. 
Set $P=<L[-n]U: n=3,4,...>$. Clearly $P\subseteq C_{[2]}(U)$ since 
$(m-1)a[-m]=(L[-1]a)[-m+1]$. $P$ is also 
invariant under the action of $Vir_-$, 
where $Vir_-=\oplus_{n=1}^{\infty}\C L[-n]$. 
Since $[L[-1], L[-2]]=L[-3]$, 
$U$ is spanned by $\{L[-1]^nL[-2]^me+ P: n,m\geq 0\}$ 
and spanned by $\{L[-2]^me+C_{[2,0]}(U): m\geq 0\}$.
 
$L(c,0)$ is  a quotient of the corresponding verma module 
$M=M(c,0)$ and 
we have $M\cong U(Vir_-)\cdot\1$ (cf. \cite{FZ}). 
We have $L(c,0)=M/J$ and $J$ contains a singular vectors of the form 
$$ \alpha=L[-2]^{m}\1+\sum a_{n_1,...,n_r}L[-n_1-2]...L[-n_r-2]\1 $$
by \cite{FF}, 
where the sum ranges over certain 
$(n_1,...,n_r)\in \Z_+^r$ with $n_1+...+n_r\not=0, 
a_{n_1,...,n_r}\in \C$. 
We note that $(L[-n]w)[-1]u\equiv L[-n]w[-1]u \pmod{P}$ for $w\in L(c,0)$, 
$u\in U$ and $n\geq 2$ by Associativity (2.1). 
Therefore, we have  
$$0=\alpha[-1]e\equiv L[-2]^me+\sum a_{n_1,...,n_r}L[-n_1-2]...L[-n_r-2]e
\equiv L[-2]^me \pmod{P}.$$  
So $U$ is spanned by $\{L[-2]^je+C_{[2,0]}(U): j=0,1,...,m-1\}$, which 
implies that $U$ satisfies condition $C_{[2,0]}$.
\prend

\begin{lmm}
Assume $U$ and $W$ are irreducible $L(c,0)$-modules and 
$I{W\choose U\quad W}\not=0$.  Let $u$ be a highest weight vector of $U$. 
Then $S^I(u,\tau)\not=0$. 
\end{lmm}

\pr
Set $U=L(c,k)$ and $W=L(c,h)$.  We note that $\dim W(h)=\dim U(k)=1$.  
Assume $S^I(u,\tau)=0$.  We first claim that $S^I(v,\tau)=0$ for 
all $v\in U$. Since $U[k]=U(k)=\C u$, $S^I(v,\tau)=0$ for $v\in U[k]$. 
Assume $S^I(w,\tau)=0$ for $[\wt]w<n+k$ and $[\wt]v=n+1+k$. 
Since $U$ is a highest weight module, we may assume 
$v=L[-1]w$ or $v=L[-2]w$ for some $w\in U$. 
Since $S^I(L[-1]w,\tau)=S^I(\tilde{\omega}[0]w,\tau)=0$ by (3.8), 
we may assume $v=L[-2]w$.  Then (3.11) implies 
$S^I(v,\tau)=(q\frac{d}{dq}-c/24)S^I(w,\tau)
+\sum_{r=1}^{\infty}E_{2r}S^I(L[2r-2]w,\tau)=0$ by 
induction. So we have $S^I(v,\tau)=0$ for $v\in U$. 
In particular, $\tr_{|W(0)} o(v)=0$ and so $o(v)_{|W(0)}=0$.   
On the other hand, by Theorem 2.10 there is a natural linear 
isomorphism 
$\pi:I{W\choose U\quad W} \to \Hom_{A(V)}(A(U)\otimes_{A(V)}W(0), W(0))$ 
satisfying $\pi(I) (v\otimes b)=o(v)b$ for $v\in A(V)$ and $b\in W(0)$. 
Therefore, 
we have $\pi(I)=0$ and so $I=0$, a contradiction. 
\prend

We will calculate some trace function explicitly. 
We always take $u\in U$ such that the coefficient 
of the leading term 
of $S(u,\tau)$ is one. 
 
For example, $L(\hf,0)$ is the first one in the discrete series and 
it has three irreducible 
modules $L(\hf,0), L(\hf,\hf), L(\hf,\st)$ as we mentioned in 
the introduction. For $U=L(\hf,\hf)$, 
$W=L(\hf,\st)$ is the only irreducible $L(\hf,0)$-module satisfying 
$0\not=I{W\choose U\quad W}$. It also satisfies 
$\dim I{W \choose U \quad W}=1$. 
Hence $S^I(v,\tau)$ is a modular form (with a linear character) 
of weight $\hf$ for $0\not=v\in L(\hf,\hf)(0)$. By the definition 
of trace function, the leading term is 
$q^{-\frac{1}{24}\frac{1}{2}+\frac{1}{16}}=q^{\frac{1}{24}}$. 
Since $S(u,\tau)$ and $\eta(\tau)$ are modular forms with 
liner characters and same leading terms and 
$S(u,\tau)/\eta(\tau)$ is holomorphic on $\CH$, 
$S(u,\tau)=\eta(\tau)$. 

The second one is $L(\frac{7}{10},0)$. 
It has 6 irreducible modules 
$L(\svt,0)$, $L(\svt,\frac{1}{10})$, 
$L(\svt,\frac{3}{5})$, $L(\svt,\frac{3}{2})$, 
$L(\svt,\frac{7}{16})$ and $L(\svt,\frac{3}{80})$. 
For $U=L(\svt,\frac{1}{10})$, $W=L(\svt,\frac{3}{80})$ is the only 
irreducible module satisfying 
$I{W\choose U\quad W}\not=0$. It also satisfies  
$\dim I{W\choose U\quad W}=1$. Hence $S^I(u,\tau)$ is a modular 
form (with a linear character) of weight $\frac{1}{10}$ for $
u\in U(0)$.  Its leading term 
$q^{-\frac{1}{24}\frac{7}{10}+\frac{3}{80}}=q^{\frac{1}{120}}$  
is equal to one of $(\eta(\tau))^{1/5}$.  Hence 
$S(u,\tau)=(\eta(\tau))^{1/5}$. 

$L(\frac{4}{5},0)$ is the third and has 10 irreducible modules.  
For $U=L(\frac{4}{5}, \frac{2}{5})$, 
$W=L(\frac{4}{5}, \frac{1}{15})$ is the only 
irreducible module satisfying 
$I{W\choose U\quad W}\not=0$. It also satisfies 
$\dim I{W\choose U\quad W}=1$. Hence $S^I(u,\tau)$ is a modular 
form (with a linear character) of weight $\frac{2}{5}$ 
for $u\in U(0)$ and its leading term is 
$q^{1\frac{1}{24}\frac{4}{5}+\frac{1}{15}}=q^{\frac{1}{30}}$, 
which is equal to the one of $(\eta(\tau))^{4/5}$. 
Hence $S(u,\tau)=(\eta(\tau))^{4/5}$.

$L(\frac{6}{7},0)$ is the fourth, which 
has 15 irreducible modules. 
For $U=L(\frac{6}{7}, \frac{1}{7})$, 
$W=L(\frac{6}{7}, \frac{1}{21})$ is the only 
irreducible module satisfying 
$I{W\choose U\quad W}\not=0$. We also have 
$\dim I{W\choose U\quad W}=1$.  
Hence $S^I(u,\tau)$ is a modular form 
(with a linear character) of weight $\frac{1}{7}$ 
for $u\in U(0)$. The leading term is 
$q^{-\frac{1}{24}\frac{6}{7}+\frac{1}{21}}=q^{1/81}$, 
which is equal to the one of $(\eta(\tau))^{2/7}$.
Hence $S(u,\tau)=(\eta(\tau))^{2/7}$.

\end{document}